\documentclass[11pt]{article}
\usepackage{amscd}
\usepackage{amsfonts}
\usepackage{amsmath}
\usepackage{amssymb}
\usepackage{amsthm}
\usepackage{bbm}
\usepackage{CJK}
\usepackage{fancyhdr}
\usepackage{graphicx}
\usepackage{hyperref}
\usepackage{indentfirst}
\usepackage{latexsym}
\usepackage{mathrsfs}
\usepackage{xypic}
\usepackage[compress]{cite}

\newtheorem{theorem}{Theorem}[section]
\newtheorem{lemma}[theorem]{Lemma}
\newtheorem{definition}[theorem]{Definition}
\newtheorem{proposition}[theorem]{Proposition}
\newtheorem{example}[theorem]{Example}

\newtheorem{cor}[theorem]{Corollary}
\newtheorem{remark}[theorem]{Remark}

\usepackage[top=1in,bottom=1in,left=1.25in,right=1.25in]{geometry}
\def\<{\langle}
\def\>{\rangle}

\def\c{\cdot}

\def\o{\otimes}

\date{}
\begin{document}
\renewcommand{\baselinestretch}{1.2}
\renewcommand{\arraystretch}{1.0}
\title{\bf Relative Rota-Baxter systems on Leibniz algebras}
\date{}
\author{{\bf Apurba Das$^{1}$,   Shuangjian Guo$^{2}$\footnote
        { Corresponding author} }\\
{\small 1. Department of Mathematics and Statistics, Indian Institute of Technology} \\
{\small  Kanpur 208016, Uttar Pradesh, India}\\
{\small Email: apurbadas348@gmail.com}\\
{\small 2. School of Mathematics and Statistics, Guizhou University of Finance and Economics} \\
{\small  Guiyang  550025, P. R. of China}\\
  {\small Email: shuangjianguo@126.com}}
 \maketitle
\begin{center}
\begin{minipage}{13.cm}

{\bf \begin{center} ABSTRACT \end{center}}
In this paper, we introduce relative Rota-Baxter systems on Leibniz algebras and give some characterizations and new constructions. Then we construct a graded Lie algebra whose Maurer-Cartan elements are relative Rota-Baxter systems. This allows us to define a cohomology theory associated with a relative Rota-Baxter system. Finally, we study formal deformations and extendibility of finite order deformations of a relative Rota-Baxter system in terms of the cohomology theory.
 \smallskip

 \smallskip

 \end{minipage}
 \end{center}
 \normalsize\vskip0.5cm

\section*{Introduction}
\def\theequation{\arabic{section}. \arabic{equation}}
\setcounter{equation} {0}

The notion of Rota-Baxter operators on associative algebras was introduced in
1960 by Baxter \cite{B60} in his study of fluctuation theory in probability. Recently, it has
been found many applications, including in Connes-Kreimer's algebraic approach
to the renormalization in perturbative quantum field theory \cite{CK00}. In the Lie algebra
context, a Rota-Baxter operator of weight zero was introduced independently in
the 1980's as the operator form of the classical Yang-Baxter equation, whereas the
classical Yang-Baxter equation plays important roles in many fields in mathematics
and mathematical physics such as quantum groups and integrable systems \cite{CP95, S83}.

\medskip

The notion of a Leibniz algebra was introduced by Bloh \cite{bloh} and rediscovered by Loday \cite{L93, LP93} in the study of the periodicity in algebraic $K$-theory. Leibniz algebras were studied from different aspects due to applications in both mathematics and physics. In particular, integration of Leibniz algebras was studied in \cite{BW17, C13} and deformation quantization of Leibniz algebras was considered in \cite{DW15}. As the underlying structure of embedding tensor, Leibniz algebras also have application in higher gauge theories, see \cite{KS20, SW20} for more details. Recently, relative Rota-Baxter operators on Leibniz algebras were studied in \cite{ST19}, which is the main ingredient in the study of the twisting theory and the bialgebra theory for Leibniz algebras. Moreover, relative Rota-Baxter operators on a Leibniz algebra can be seen as the Leibniz algebraic analogue of Poisson structures. Generally, Rota-Baxter operators can be defined on operads, which give rise to the splitting of operads \cite{BBGN13, PBG17}. For further details on Rota-Baxter operators, see \cite{G12}.

\medskip

The deformation of algebraic structures began with the seminal work of Gerstenhaber \cite{G63, G64} for associative algebras and followed by its extension to Lie algebras by Nijenhuis and Richardson \cite{NR66, NR68}.  In general, the deformation theory of algebras over binary quadratic operads was developed by Balavoine \cite{B97}. Deformations of morphisms and $\mathcal{O}$-operators (also called relative Rota-Baxter operators) were developed in \cite{Da20, FZ15} and \cite{TBGS19, TBGS20}. Rota-Baxter systems as a generalization of a Rota-Baxter operator were introduced by Brzezi\'{n}ski \cite{B16}. In a Rota-Baxter system, two operators are acting on the algebra and satisfying some Rota-Baxter type identities. Rota-Baxter systems in the presence of bimodule were introduced and their deformation theory was studied by Das \cite{D20}. They are called generalized Rota-Baxter systems. Our main objective in this paper is generalized Rota-Baxter systems in the context of Leibniz algebra. We call them relative Rota-Baxter systems, motivated by the terminology of relative Rota-Baxter operators of \cite{ST19}. Our aim in this paper is to study the cohomology and deformation theory of relative Rota-Baxter systems in the context of Leibniz algebras.

\medskip

The paper is organized as follows. In Section 2, we introduce relative Rota-Baxter systems with respect to a representation of a Leibniz algebra and give some characterizations and new constructions. In Section 3, we emphasis on relative Rota-Baxter systems with respect to the regular representation. In Section 4, we construct a graded Lie algebra whose Maurer-Cartan elements are relative Rota-Baxter systems, which leads us to define cohomology for a relative Rota-Baxter system. Finally, in Section 5, we consider formal deformations of relative Rota-Baxter systems.

\medskip

Throughout this paper, $\mathbb{K}$ is a field of characteristic zero and $\mathbb{Z}$ denotes the set of all integers.

\section{Leibniz algebras}
\def\theequation{\arabic{section}.\arabic{equation}}
\setcounter{equation} {0}

In this preliminary section, we recall Leibniz algebras and their representations \cite{L93, LP93}.

\begin{definition}
A Leibniz algebra is a vector space $g$ together with a bilinear operation $[\c, \c]_g : g \o g \rightarrow  g$ satisfying
\begin{eqnarray*}
[x, [y, z]_g]_g=[[x, y]_g, z]_g+[y, [x, z]_g]_g, ~ \text{ for } x, y, z\in g.
\end{eqnarray*}
\end{definition}

\begin{definition}
A representation of a Leibniz algebra $(g, [\c, \c]_g)$ is a triple $(V,  \rho^L, \rho^R)$, where $V$ is a vector space,
$\rho^L, \rho^R: g \rightarrow gl(V)$ are linear maps such that the following equalities hold :  for all $x, y \in g$,
\begin{eqnarray*}
(1)~~&&\rho^L([x, y]_g) = [\rho^L(x), \rho^L(y)],\\
(2)~~&&\rho^R([x, y]_g) = [\rho^L(x), \rho^R(y)], \\
 (3)~~&&\rho^R(y) \circ \rho^L(x) = -\rho^R(y) \circ \rho^R(x).
\end{eqnarray*}
\end{definition}

Let $(g, [\c, \c])$ be a Leibniz algebra.
Define the left multiplication $L: g \rightarrow gl(g)$ and the right multiplication $R: g \rightarrow gl(g)$ by
$L_xy = [x, y]_g$ and $R_xy = [y, x]_g$, for all $x, y \in g$. Then $(g, L, R)$ is a representation of
$(g, [\c, \c]_g)$, called the regular representation. Define two linear maps $L^{\ast}, R^{\ast}: g \rightarrow gl(g^{\ast})$
with $x \mapsto L^{\ast}_x$ and $x \mapsto R^{\ast}_x$ respectively by
\begin{eqnarray*}
\langle L^{\ast}_x\xi, y\rangle =-\langle \xi, [x, y]_g\rangle,~~~\langle R^{\ast}_x\xi, y\rangle =-\langle \xi, [y, x]_g\rangle, ~ \text{ for } x,y\in g,~ \xi\in g^{\ast}.
\end{eqnarray*}

Then it has been shown in \cite{ST19} that $(g^*, L^*, -L^* - R^*)$ is a representation. This is called the dual of the regular representation.

\begin{definition}
A quadratic Leibniz algebra is a Leibniz algebra $(g, [\c, \c]_g)$ equipped with a nondegenerate skew-symmetric bilinear form $\omega\in \wedge^2g^{\ast}$ such that the following invariant condition holds:
\begin{eqnarray*}
\omega(x, [y, z]_g)=\omega([x, z]_g+[z,x]_g, y),~ \text{ for } x, y, z \in g.
\end{eqnarray*}
\end{definition}

\begin{proposition} (\cite{ST19})
Let $(g, [\c, \c]_g, \omega)$ be a quadratic Leibniz algebra. Then the map $$\omega^{\natural}:g\rightarrow g^{\ast}, ~ \omega^\sharp (x) (y) = \omega (x, y), ~ \text{ for } x, y \in g$$
is an isomorphism from the regular representation $(g,  L, R)$ to its dual representation $(g^{\ast}, L^{\ast}, - L^{\ast}-R^{\ast})$.
\end{proposition}

\section{Relative Rota-Baxter systems with respect to a representation}
\def\theequation{\arabic{section}. \arabic{equation}}
\setcounter{equation} {0}

In this section, we introduce relative Rota-Baxter systems with respect to a representation of a Leibniz algebra.

Let $(V, \rho^L, \rho^R)$ be a representation of a Leibniz algebra $(g, [\c, \c]_g)$.

\begin{definition}
(1) A relative Rota-Baxter system on  $(g, [\c, \c]_g)$ with respect to the representation $(V, \rho^L, \rho^R)$ consists of a pair $(R,S)$ of linear maps  $R, S : V \rightarrow g$
satisfying
\begin{align*}
[Ru, Rv]_g =~ R(\rho^{L}(Ru)v+\rho^{R}(Sv)u),\\
[Su, Sv]_g =~ S(\rho^{L}(Ru)v+\rho^{R}(Sv)u),
\end{align*}
for $u, v \in V$.

\medskip

(2) A Rota-Baxter system on $(g, [\c, \c]_g)$ is a relative Rota-Baxter system on $(g, [\c, \c]_g)$ with respect to the regular representation.
\end{definition}

\begin{example}
A relative Rota-Baxter operator \cite{ST19} on $(g, [\c, \c]_g)$ with respect to the representation $(V, \rho^L, \rho^R)$ is a linear map $R : V \rightarrow g$ satisfying
\begin{align*}
    [Ru, Rv]_g = R(\rho^{L}(Ru)v+\rho^{R}(Rv)u), ~\text{ for } u, v \in V.
\end{align*}
Thus $R$ is a relative Rota-Baxter operator if and only if the pair $(R,R)$ is a relative Rota-Baxter system.
\end{example}

\begin{example}
 Consider the 2-dimensional Leibniz algebra $(g, [\c, \c])$ given with respect to a basis
$\{e_1, e_2\}$ by
\begin{eqnarray*}
[e_1,e_1]=0,~~~ [e_1,e_2]=0,~~~ [e_2,e_1]=e_1,~~~[e_2,e_2]=e_1.
\end{eqnarray*}
Let $\{e^{\ast}_1, e^{\ast}_2\}$ be the dual basis. Then $R=\left(
                                                              \begin{array}{cc}
                                                                a_{11} & a_{12} \\
                                                                a_{21} & a_{22} \\
                                                              \end{array}
                                                            \right), ~S=\left(
                                                              \begin{array}{cc}
                                                                b_{11} & b_{12} \\
                                                                b_{21} & b_{22} \\
                                                              \end{array}
                                                            \right)$ is a relative Rota-Baxter system on $(g, [\c, \c])$
with respect to the representation $(g^{\ast}, L^{\ast}, - L^{\ast}-R^{\ast})$ if and only if
\begin{eqnarray*}
&&[Re_i^{\ast},Re_j^{\ast} ]=R(L^{\ast}_{Re_i^{\ast}}e_j^{\ast}-L^{\ast}_{Se_j^{\ast}}e_i^{\ast}-R^{\ast}_{Se_j^{\ast}}e_i^{\ast}),\\
&& [Se_i^{\ast},Se_j^{\ast} ]=S(L^{\ast}_{Re_i^{\ast}}e_j^{\ast}-L^{\ast}_{Se_j^{\ast}}e_i^{\ast}-R^{\ast}_{Se_j^{\ast}}e_i^{\ast}),~~~i,j=1, 2.
\end{eqnarray*}
It is straightforward to deduce that
\begin{eqnarray*}
&& L_{e_1}(e_1,e_2)=(e_1,e_2)\left(
                                              \begin{array}{cc}
                                                                0 & 0 \\
                                                                0 & 0 \\
                                                              \end{array}
                                                            \right),~~~ L_{e_2}(e_1,e_2)=(e_1,e_2)\left(
                                              \begin{array}{cc}
                                                                1 & 1 \\
                                                                0 & 0 \\
                                                              \end{array}
                                                            \right),\\
&&  R_{e_1}(e_1,e_2)=(e_1,e_2)\left(
                                              \begin{array}{cc}
                                                                0 & 1 \\
                                                                0 & 0 \\
                                                              \end{array}
                                                            \right),~~~ R_{e_2}(e_1,e_2)=(e_1,e_2)\left(
                                              \begin{array}{cc}
                                                                0 & 1 \\
                                                                0 & 0 \\
                                                              \end{array}
                                                            \right),
\end{eqnarray*}
and
\begin{eqnarray*}
&& L^{\ast}_{e_1}(e^{\ast}_1,e^{\ast}_2)=(e^{\ast}_1,e^{\ast}_2)\left(
                                              \begin{array}{cc}
                                                                0 & 0 \\
                                                                0 & 0 \\
                                                              \end{array}
                                                            \right),~~~~~ L^{\ast}_{e_2}(e^{\ast}_1,e^{\ast}_2)=(e^{\ast}_1,e^{\ast}_2)\left(
                                              \begin{array}{cc}
                                                                -1 & 0 \\
                                                                -1 & 0 \\
                                                              \end{array}
                                                            \right),\\
&&  R^{\ast}_{e_1}(e^{\ast}_1,e^{\ast}_2)=(e^{\ast}_1,e^{\ast}_2)\left(
                                              \begin{array}{cc}
                                                                0 & 0 \\
                                                                -1 & 0 \\
                                                              \end{array}
                                                            \right),~~~ R_{e_2}^{\ast}(e^{\ast}_1,e^{\ast}_2)=(e^{\ast}_1,e^{\ast}_2)\left(
                                              \begin{array}{cc}
                                                                0 & 0 \\
                                                                -1 & 0 \\
                                                              \end{array}
                                                            \right).
\end{eqnarray*}
We have
\begin{eqnarray*}
[Re^{\ast}_1,Re^{\ast}_1]=[a_{11}e_1+a_{21}e_2, a_{11}e_1+a_{21}e_2]=a_{21}(a_{11}+a_{21})e_1,
\end{eqnarray*}
and
\begin{eqnarray*}
&&R(L^{\ast}_{Re_1^{\ast}}e_1^{\ast}-L^{\ast}_{Se_1^{\ast}}e_1^{\ast}-R^{\ast}_{Se_1^{\ast}}e_1^{\ast})\\
&=&- a_{21}(R(e_1^{\ast})+R(e_2^{\ast})) +  b_{21}(R(e_1^{\ast})+R(e_2^{\ast})) + (b_{11} + b_{21})R(e_2^{\ast}) \\
&=&- a_{21}(a_{11}e_1+a_{21}e_2+a_{12}e_1+a_{22}e_2)+  b_{21}(a_{11}e_1+a_{21}e_2+a_{12}e_1+a_{22}e_2)\\
&& +(b_{11}+ b_{21})(a_{12}e_1+a_{22}e_2)\\
&=& ((b_{11} + b_{21})a_{12}+(a_{11} + a_{12})(b_{21}-a_{21}))e_1+((b_{11} + b_{21})a_{22}+(a_{21}+a_{22})(b_{21}-a_{21}))e_2,
\end{eqnarray*}
\begin{eqnarray*}
[Se^{\ast}_1,Se^{\ast}_1]=[b_{11}e_1+b_{21}e_2, b_{11}e_1+b_{21}e_2]=b_{21}(b_{11}+b_{21})e_1,
\end{eqnarray*}
and
\begin{eqnarray*}
&&S(L^{\ast}_{Re_1^{\ast}}e_1^{\ast}-L^{\ast}_{Se_1^{\ast}}e_1^{\ast}-R^{\ast}_{Se_1^{\ast}}e_1^{\ast})\\
&=& - a_{21}(R(e_1^{\ast})+R(e_2^{\ast})) + b_{21}(S(e_1^{\ast})+S(e_2^{\ast})) + (b_{11} + b_{21})S(e_2^{\ast}) \\
&=&- a_{21}(a_{11}e_1+a_{21}e_2+a_{12}e_1+a_{22}e_2)+  b_{21}(b_{11}e_1+b_{21}e_2+b_{12}e_1+b_{22}e_2)\\
&& +(b_{11}+ b_{21})(b_{12}e_1+b_{22}e_2)\\
&=& ((b_{11} + b_{21})b_{12}+(b_{11} + b_{12})b_{21}-(a_{11} + a_{12})a_{21})e_1\\
&&+((b_{11} + b_{21})b_{22}+(b_{21} + b_{22})b_{21}-(a_{21} + a_{22})a_{21})e_2.
\end{eqnarray*}
Thus, we obtain
\begin{eqnarray*}
&&a_{21}(a_{11}+a_{21})=(b_{11} + b_{21})a_{12}+(a_{11} + a_{12})(b_{21}-a_{21}), \\
&& (b_{11} + b_{21})a_{22}+(a_{21}+a_{22})(b_{21}-a_{21})=0,\\
&&b_{21}(b_{11}+b_{21})=(b_{11} + b_{21})b_{12}+(b_{11} + b_{12})b_{21}-(a_{11} + a_{12})a_{21},\\
&& (b_{11} + b_{21})b_{22}+(b_{21} + b_{22})b_{21}-(a_{21} + a_{22})a_{21}=0.
\end{eqnarray*}
Similarly, we obtain
\begin{eqnarray*}
&& a_{21}(a_{12}+a_{22})=b_{22}(a_{11}+a_{12})+(b_{12}+b_{22})a_{12},\\
&& b_{22}(a_{21}+a_{22})+(b_{12}+b_{22})a_{22}=0,\\
&& b_{21}(b_{12}+b_{22})=b_{22}(b_{11}+b_{12})+(b_{12}+b_{22})b_{12},\\
&&b_{22}(b_{21}+b_{22})+(b_{12}+b_{22})b_{22}=0,~~~~ a_{22}(a_{12}+a_{22})=b_{22}(b_{12}+b_{22})=0,\\
&& a_{22}(a_{11} + a_{21}) =a_{22}(a_{11} + a_{12}),~~~~ -a_{22}(a_{21} + a_{22}) = 0,\\
&& b_{22}(b_{11} + b_{21}) =a_{22}(b_{11} + b_{12}),~~~~ -a_{22}(b_{21} + b_{22}) = 0.
\end{eqnarray*}
Summarize the above discussion, we have

(1) If $a_{22}=b_{22}=0$ and $a_{21}=b_{21}$, then $R=\left(
                                                              \begin{array}{cc}
                                                                a_{11} & a_{12} \\
                                                                a_{21} & 0 \\
                                                              \end{array}
                                                            \right), S=\left(
                                                              \begin{array}{cc}
                                                                b_{11} & b_{12} \\
                                                                b_{21} & 0 \\
                                                              \end{array}
                                                            \right)$ is a relative Rota-Baxter system on $(g, [\c, \c])$
with respect to the representation $(g^{\ast}, L^{\ast}, -L^{\ast}-R^{\ast})$ if and only if
\begin{eqnarray*}
&&(b_{12}-a_{21})a_{12} =(b_{12}-a_{21})b_{12}=0,\\
&&a_{21}(a_{11}+a_{21})=(b_{11} + b_{21})a_{12}, \\
&&b_{21}(b_{11}+b_{21})=(b_{11} + b_{21})b_{12}+(b_{11} + b_{12})b_{21}-(a_{11} + a_{12})a_{21},
\end{eqnarray*}

(2) If $a_{22}=b_{22}\neq0$ and $a_{21}\neq b_{21}$, then $R=\left(
                                                              \begin{array}{cc}
                                                                a_{11} & a_{12} \\
                                                                a_{21} & a_{22} \\
                                                              \end{array}
                                                            \right), S=\left(
                                                              \begin{array}{cc}
                                                                b_{11} & b_{12} \\
                                                                b_{21} & b_{22}  \\
                                                              \end{array}
                                                            \right)$ is a relative Rota-Baxter system on $(g, [\c, \c])$
with respect to the representation $(g^{\ast}, L^{\ast}, -L^{\ast}-R^{\ast})$ if and only if
\begin{eqnarray*}
a_{11}=-a_{12}=-a_{21} = a_{22},  ~~~~~~b_{11}=-b_{12}=-b_{21} = b_{22}.
\end{eqnarray*}
\end{example}

We will give some more examples of Rota-Baxter systems on Leibniz algebras in the next section.

\medskip

In the following, we give some characterizations of relative Rota-Baxter systems.
Let $(V, \rho^L, \rho^R)$ be a representation of a Leibniz algebra $(g, [\c, \c]_g)$.  Then there is a Leibniz algebra structure  on $g\oplus g\oplus V$  given by
\begin{eqnarray*}
[x_1+x_2+u, y_1+y_2+v]=[x_1, y_1]_g+[x_2, y_2]_g+\rho^{L}(x_1)v+\rho^{R}(y_2)u.
\end{eqnarray*}
This is exactly the semidirect product if we consider the Leibniz algebra structure on   $g\oplus g$ and define its representation on $V$ by
$\rho^{L} (x_1+x_2)v = \rho^{L} (x_1)v$ and $\rho^{R} (x_1+x_2)v = \rho^{R} (x_2)v$.

\begin{proposition}
A pair $(R, S)$ of linear maps from $V$ to $g$ is a relative Rota-Baxter system with respect to the representation $(V, \rho^L, \rho^R)$  if and only if the pair $(\widetilde{R}, \widetilde{S})$ of maps
\begin{eqnarray*}
&& \widetilde{R}: g\oplus g\oplus V\rightarrow  g\oplus g\oplus V, ~~~x_1+x_2+u\mapsto R(u)+0+0,\\
&& \widetilde{S}: g\oplus g\oplus V\rightarrow  g\oplus g\oplus V,~~~x_1+x_2+u\mapsto 0+S(u)+0,
\end{eqnarray*}
is a Rota-Baxter system on the Leibniz algebra $g\oplus g\oplus V$.
\end{proposition}
{\bf Proof.}  For any $x_1, x_2, y_1, y_2\in g$ and $u, v\in V$,  we have
\begin{eqnarray*}
[\widetilde{R}(x_1+x_2+u), \widetilde{R}(y_1+y_2+v)]=[R(u), R(v)]_g + 0 + 0
\end{eqnarray*}
and
\begin{eqnarray*}
&&\widetilde{R}([\widetilde{R}(x_1+x_2+u), y_1+y_2+v]+[x_1+x_2+u, \widetilde{S}(y_1+y_2+v)])\\
&=& R(\rho^{L}(Ru)v+\rho^{R}(Sv)u) + 0 + 0.
\end{eqnarray*}
Similarly,  we have \begin{eqnarray*}
[\widetilde{S}(x_1+x_2+u), \widetilde{S}(y_1+y_2+v)]= 0 + [S(u), S(v)]_g + 0
\end{eqnarray*}
and
\begin{eqnarray*}
&&\widetilde{S}([\widetilde{R}(x_1+x_2+u), y_1+y_2+v]+[x_1+x_2+u, \widetilde{S}(y_1+y_2+v)]\\
&=& 0 + S(\rho^{L}(Ru)v+\rho^{R}(Sv)u) + 0.
\end{eqnarray*}
Hence $(R, S)$ is a relative Rota-Baxter system if and only if $( \widetilde{R},  \widetilde{S})$ is a Rota-Baxter system. \hfill $\square$

Recall that a Nijenhuis operator
on a Leibniz algebra $(g, [\c, \c]_g)$ is a linear map $N : g\rightarrow g$ satisfying
\begin{eqnarray*}
[Nx, Ny]_g=N([N(x), y]_g + [x, N(y)]_g-N[x, y]_g ),~ \text{ for } x, y\in g.
\end{eqnarray*}
The following result relates to relative Rota-Baxter systems and Nijenhuis operators.

\begin{proposition}
A pair $(R, S)$ of linear maps from $V$ to $\mathfrak{g}$ is a relative Rota-Baxter system if and only if
\begin{eqnarray*}
N_{(R, S)}=\left(
                                            \begin{array}{ccc}
                                              0 & 0 & R \\
                                             0 & 0 & S \\
                                              0 & 0 & 0 \\
                                            \end{array}
                                          \right): g\oplus g\oplus V\rightarrow  g\oplus g\oplus V
\end{eqnarray*}
is a Nijenhuis operator on the Leibniz algebra $g\oplus g\oplus V$.
\end{proposition}
{\bf Proof.} For any $x_1, x_2, y_1, y_2\in g$ and $u, v\in V$, by a simple calculation, we have
\begin{eqnarray*}
[N_{(R, S)}(x_1+y_1+u), N_{(R, S)}(x_2+y_2+v)] =[R(u), R(v)]_g+ [S(u), S(v)]_g + 0
\end{eqnarray*}
and
\begin{eqnarray*}
&&N_{(R, S)}([N_{(R, S)}(x_1+y_1+u), x_2+y_2+v] + [x_1+y_1+u, N_{(R, S)}(x_2+y_2+v)]\\
&&-N_{(R, S)}[x_1+y_1+u, x_2+y_2+v] )\\
&=& R(\rho^{L}(Ru)v+\rho^{R}(Sv)u)+S(\rho^{L}(Ru)v+\rho^{R}(Sv)u) + 0.
\end{eqnarray*}
It follows  that $N_{(R, S)}$ is a Nijenhuis operator if and only if $(R, S)$ is a relative Rota-Baxter system. \hfill $\square$

\begin{definition}
Let $(V, \rho^L, \rho^R)$ be a representation of a Leibniz algebra $(g, [\c, \c]_g)$.  Suppose that  $dim ~(g) = dim~ (V)$. A pair $(\Phi, \Psi)$
of invertible linear maps from $g$ to $V$ is said to be an invertible 1-cocycle system if they satisfy
\begin{eqnarray*}
&&\Phi([x, y]_g)=\rho^{L}(x)\Phi(y)+\rho^{R}(\Psi^{-1}\circ \Phi(y))\Phi(x),\\
&& \Psi([x, y]_g)=\rho^{L}(\Phi^{-1}\circ \Psi(x))\Psi(y)+\rho^{R}(y)\Psi(x),
\end{eqnarray*}
for $x, y \in g$.
\end{definition}
 It follows from the above definition that $(\Phi, \Phi)$ is an invertible 1-cocycle system if and only
if $\Phi : g \rightarrow V$ is an invertible derivation.

\begin{proposition}
Let $(V, \rho^L, \rho^R)$ be a representation of a Leibniz algebra $(g, [\c, \c]_g)$. Suppose that $dim~(g) = dim~( V)$.  A pair $(R, S)$
of invertible linear maps from $V$ to $g$ is a relative Rota-Baxter system if and only if $(R^{-1}, S^{-1})$ is an invertible 1-cocycle system.
\end{proposition}
{\bf Proof.} For any $u, v\in V$ and $x, y\in g$, by taking $R(u) = x, R(v) = y$, the first identity of Definition 2.1 is equivalent to
\begin{eqnarray*}
&&R^{-1}[x, y]_g=\rho^{L}(x)R^{-1}{y}+\rho^{R}((S^{-1})^{-1}\circ R^{-1}(y))R^{-1}{x}.
\end{eqnarray*}
Similarly, for any $u, v\in V$ and $x, y\in g$, by taking $S(u) = x, ~S(v) = y$, the second identity of Definition 2.1 is equivalent to
\begin{eqnarray*}
&&S^{-1}[x, y]_g=\rho^{L}((R^{-1})^{-1}\circ S^{-1}(x))S^{-1}{y}+\rho^{R}(y)S^{-1}{x}.
\end{eqnarray*}
It follows that $(R, S)$
of invertible linear maps from $V$ to $g$ is a relative Rota-Baxter system if and only if $(R^{-1}, S^{-1})$ is an invertible 1-cocycle system. \hfill $\square$

\medskip

\begin{proposition}
Let $(g, [\c, \c]_g, \omega)$ be a quadratic Leibniz algebra and $R, S: g^{\ast}\rightarrow g$ be two linear
maps. Then $(R, S)$ is a relative Rota-Baxter system on $(g, [\c, \c]_g)$ with respect to the representation
$(g^{\ast}, L^{\ast}, - L^{\ast}-R^{\ast})$ if and only if $(R\circ\omega^{\natural}, S\circ \omega^{\natural})$ is a Rota-Baxter system on $(g, [\c, \c]_g)$.
\end{proposition}

{\bf Proof.} For any $x, y\in g$, we have
\begin{eqnarray*}
&&R\circ \omega^{\natural}([R\circ\omega^{\natural}(x), y]_g+[x, S\circ\omega^{\natural}(y)]_g)\\
&=&R(\omega^{\natural}(L_{R\circ\omega^{\natural}(x)}y)+\omega^{\natural}(R_{S\circ\omega^{\natural}(y)}x))\\
&=&R(L^{\ast}_{R\circ\omega^{\natural}(x)}\omega^{\natural}(y)-L^{\ast}_{S\circ\omega^{\natural}(y)}\omega^{\natural}(x)-R^{\ast}_{S\circ\omega^{\natural}(y)}\omega^{\natural}(x)).
\end{eqnarray*}
Similarly, we have
\begin{eqnarray*}
&&S\circ \omega^{\natural}([R\circ\omega^{\natural}(x), y]_g+[x, S\circ\omega^{\natural}(y)]_g)\\
&=&S(\omega^{\natural}(L_{R\circ\omega^{\natural}(x)}y)+\omega^{\natural}(R_{S\circ\omega^{\natural}(y)}x))\\
&=&S(L^{\ast}_{R\circ\omega^{\natural}(x)}\omega^{\natural}(y)-L^{\ast}_{S\circ\omega^{\natural}(y)}\omega^{\natural}(x)-R^{\ast}_{S\circ\omega^{\natural}(y)}\omega^{\natural}(x)).
\end{eqnarray*}
Thus it follows that $(R\circ\omega^{\natural}, S\circ \omega^{\natural})$ is a Rota-Baxter system on $(g, [\c, \c]_g)$ if and only if
\begin{eqnarray*}
&&[R\circ \omega^{\natural}(x), R\circ \omega^{\natural}(y)]_g=R(L^{\ast}_{R\circ \omega^{\natural}(x)}\omega^{\natural}(y)-L^{\ast}_{S\circ \omega^{\natural}(y)}\omega^{\natural}(x)-R^{\ast}_{S\circ \omega^{\natural}(y)}\omega^{\natural}(x)),\\
&&[S\circ \omega^{\natural}(x), S\circ \omega^{\natural}(y)]_g=S(L^{\ast}_{R\circ \omega^{\natural}(x)}\omega^{\natural}(y)-L^{\ast}_{S\circ \omega^{\natural}(y)}\omega^{\natural}(x)-R^{\ast}_{S\circ\omega^{\natural}(y)}\omega^{\natural}(x)).
\end{eqnarray*}
Since $\omega^{\natural}$ is an isomorphism,  these identities hold if and only if  $(R, S)$ is a relative Rota-Baxter system on $(g, [\c, \c]_g)$ with respect to the representation
$(g^{\ast}, L^{\ast}, -L^{\ast}-R^{\ast})$. \hfill $\square$

\section{Rota-Baxter systems}
\def\theequation{\arabic{section}. \arabic{equation}}
\setcounter{equation} {0}

In this section, we mainly provide examples of Rota-Baxter systems on Leibniz algebras. As mentioned earlier, they are relative Rota-Baxter operators with respect to the regular representation.


\begin{example}
Consider the three-dimensional Leibniz algebra $(g, [\c, \c]_g)$ given with
respect to a basis $\{e_1, e_2, e_3\}$ by
\begin{eqnarray*}
[e_1, e_1]_g = e_3.
\end{eqnarray*}
Then $R=\left(
                                            \begin{array}{ccc}
                                              a_{11} & a_{12} & a_{13} \\
                                              a_{21} & a_{22} & a_{23} \\
                                              a_{31} & a_{32} & a_{33} \\
                                            \end{array}
                                          \right)$, S=$\left(
                                            \begin{array}{ccc}
                                              b_{11} & b_{12} & b_{13} \\
                                              b_{21} & b_{22} & b_{23} \\
                                              b_{31} & b_{32} & b_{33} \\
                                            \end{array}
                                          \right)$              is a Rota-Baxter system on $(g, [\c, \c]_g)$ if and only if
\begin{eqnarray*}
&&[Re_i,Re_j]_g = R ([R e_i, e_j ]_g+[e_i, Se_j]_g),  \\
&&[Se_i,Se_j]_g = S ([R e_i, e_j ]_g+[e_i, Se_j]_g),  ~ \text{ for } i, j=1, 2, 3.
\end{eqnarray*}
We have $[R e_1,Re_1]_g=[a_{11}e_1 + a_{21}e_2 + a_{31}e_3, a_{11}e_1 + a_{21}e_2 + a_{31}e_3] = a^2_{11}e_3$, and
\begin{eqnarray*}
&&R ([R e_1, e_1 ]_g+[e_1, Se_1]_g)\\
&=&R ([a_{11}e_1 + a_{21}e_2 + a_{31}e_3, e_1]_g+[e_1, b_{11}e_1 + b_{21}e_2 + b_{31}e_3]_g)\\
&=& (a_{11}+b_{11})Re_3\\
&=& (a_{11}+b_{11})a_{13}e_1+(a_{11}+b_{11})a_{23}e_2+(a_{11}+b_{11})a_{33}e_3.
\end{eqnarray*}
Thus, by $[R e_1,Re_1]_g=R ([R e_1, e_1 ]_g+[e_1, Se_1]_g)$, we have
\begin{eqnarray*}
 (a_{11}+b_{11})a_{13}=0,~~ (a_{11}+b_{11})a_{23}=0,~~ a^2_{11}=(a_{11}+b_{11})a_{33}.
\end{eqnarray*}
Similarly, by $[S e_1,Se_1]_g=S ([R e_1, e_1 ]_g+[e_1, Se_1]_g)$, we have
\begin{eqnarray*}
(a_{11}+b_{11})b_{13}=0,~~ (a_{11}+b_{11})b_{23}=0,~~ b^2_{11}=(a_{11}+b_{11})b_{33}.
\end{eqnarray*}

By considering other choices of $e_i$ and $e_j$, we obtain
\begin{eqnarray*}
&& a_{11}a_{12}=b_{12}a_{33},~~~~~~~~~~~ b_{12}a_{13}=0,~~~~~~~~~~ b_{12}a_{23}=0,\\
&& b_{11}b_{12}=b_{12}b_{33},~~~~~~~~~~~~ b_{12}b_{13}=0,~~~~~~~~~~ b_{12}b_{23}=0,\\
&& a_{11}a_{13} = b_{13}a_{33},~~~~~~~~~~~ b_{13}a_{13} = 0,~~~~~~~~~~ b_{13}a_{23} = 0,\\
&& b_{11}b_{13} = b_{13}b_{33}, ~~~~~~~~~~~~b_{13}b_{13} = 0, ~~~~~~~~~~b_{13}b_{23} = 0,\\
&& a_{12}a_{11} = a_{12}a_{33}, ~~~~~~~~~~~a_{12}a_{13} = 0, ~~~~~~~~~~a_{12}a_{23} = 0,\\
&& b_{12}b_{11} = a_{12}b_{33}, ~~~~~~~~~~~~a_{12}b_{13} = 0, ~~~~~~~~~~a_{12}b_{23} = 0,\\
&& a_{13}a_{11} = a_{13}a_{33},~~~~~~~~~~~a_{13}a_{13} = 0,~~~~~~~~ ~~a_{13}a_{23} = 0,\\
&&b_{13}b_{11} = a_{13}b_{33},~~~ ~~~~~~~~~a_{13}b_{13} = 0,~~~~~~~~~~a_{13}b_{23} = 0,\\
&& a^2_{12} = 0, ~~~~a^{2}_{13} = 0,~~~ a_{12}a_{13} = 0,~~~~ b^2_{12} = 0,~~~~ b^{2}_{13} = 0,~~~~ b_{12}b_{13} = 0.
\end{eqnarray*}
Summarize the above discussion, we have

(1) If $a_{11}=b_{11}=a_{12}=b_{12}=a_{13}=b_{13}=0$,  then any $R=\left(
                                            \begin{array}{ccc}
                                              0 & 0 & 0 \\
                                              a_{21} & a_{22} & a_{23} \\
                                              a_{31} & a_{32} & a_{33} \\
                                            \end{array}
                                          \right)$, S=$\left(
                                            \begin{array}{ccc}
                                              0 & 0 & 0\\
                                              b_{21} & b_{22} & b_{23} \\
                                              b_{31} & b_{32} & b_{33} \\
                                            \end{array}
                                          \right)$              is a Rota-Baxter system on $(g, [\c, \c]_g)$  with respect to the regular representation.

(2) If $a_{12}=b_{12}=a_{13}=b_{13}=0$ and $a_{11}=b_{11} \neq 0$, $a_{23}=b_{23}= 0$,  then any any $R=\left(
                                            \begin{array}{ccc}
                                              a_{11} & 0 & 0 \\
                                              a_{21} & a_{22} & a_{23} \\
                                              a_{31} & a_{32} & \frac{a_{11}}{2} \\
                                            \end{array}
                                          \right)$, S=$\left(
                                            \begin{array}{ccc}
                                              b_{11} & 0 & 0\\
                                              b_{21} & b_{22} & b_{23} \\
                                              b_{31} & b_{32} & \frac{b_{11}}{2} \\
                                            \end{array}
                                          \right)$              is a Rota-Baxter system on $(g, [\c, \c]_g)$  with respect to the regular representation.
\end{example}

\medskip

We have seen that relative Rota-Baxter systems generalize relative Rota-Baxter operators. In the following, we show that they also generalize Rota-Baxter operators of arbitrary weight.
\begin{definition}
Let $(g, [\c, \c]_g)$ be a Leibniz algebra. A linear map $R: g\rightarrow g$ is said to be a  Rota-Baxter operator of  weight $\lambda$ if $R$ satisfies
\begin{eqnarray*}
[Rx, Ry]_g=R([Rx, y]_g+[x, Ry]_g+\lambda[x, y]_g), ~ \text{ for } x, y\in g.
\end{eqnarray*}
\end{definition}

\begin{proposition}
Let $(g, [\c, \c]_g)$ be a Leibniz algebra and $R: g\rightarrow g$ be a   Rota-Baxter operator of  weight $\lambda$. Then $(R, R+\lambda Id)$ and $(R+\lambda Id, R)$ are Rota-Baxter systems on $(g, [\c, \c]_g)$.
\end{proposition}
{\bf Proof.} For any $x, y\in g$, we have
\begin{eqnarray*}
[Rx, Ry]_g&=& R([Rx, y]_g+[x, Ry]_g+\lambda[x, y]_g)\\
&=& R([Rx, y]_g+[x, (R+\lambda Id)y]_g) \\
&=& R ( [ (R + \lambda Id)x, y]_g + [x, Ry]_g ),
\end{eqnarray*}
and
\begin{eqnarray*}
&&[(R+\lambda Id)x, (R+\lambda Id)y]_g\\
&=& [Rx, Ry]_g+\lambda[Rx, y]_g+\lambda[x, Ry]_g+ [\lambda x, \lambda y]_g\\
&=& R([Rx, y]_g+[x, Ry]_g+\lambda[x, y]_g)+\lambda[Rx, y]_g+\lambda[x, Ry]_g+  [\lambda x, \lambda y]_g\\
&=& (R+\lambda Id)([Rx, y]_g+[x, (R+\lambda Id)y]_g)\\
&=&  (R+\lambda Id)( [ (R + \lambda Id)x, y]_g + [x, Ry]_g ).
\end{eqnarray*}
This shows that $(R, R+\lambda Id)$ and $(R+\lambda Id, R)$ are Rota-Baxter systems on $(g, [\c, \c]_g)$.
\hfill $\square$

\medskip

Let $(g, [\c, \c]_g)$ be a Leibniz algebra. A linear map $T: g\rightarrow g$ is said to be a left $g$-linear map (resp. right $g$-linear map) if $T[x, y]=[x, Ty]$ (resp. $T[x, y]=[Tx, y]$), for any $x, y \in g$.

\begin{lemma}
Let $(g, [\c, \c]_g)$ be a Leibniz algebra. Suppose that $R: g\rightarrow g$ is a  left $g$-linear map and $S: g\rightarrow g$ is a  right $g$-linear map. Then $(R, S)$ is a  Rota-Baxter system on $(g, [\c, \c]_g)$ if and only if \begin{eqnarray*}
[x, R\circ S(y)]_g=0=[S\circ R(x), y]_g,~ \text{ for } x, y\in g.
\end{eqnarray*}
\end{lemma}
{\bf Proof.} For any $x, y\in g$, we observe that
\begin{eqnarray*}
 R([Rx, y]_g+[x, Sy]_g) = [Rx, Ry]_g+[x, R\circ  S(y)]_g,
\end{eqnarray*}
and
\begin{eqnarray*}
S([Rx, y]_g+[x, Sy]_g) = [R\circ S (x), Ry]_g+[Sx, Sy]_g.
\end{eqnarray*}
It follows from the above two identities that $(R, S)$ is a Rota-Baxter system on $(g, [\c, \c]_g)$ if and only if \begin{eqnarray*}
[x, R\circ S(y)]_g=0=[S\circ R(x), y]_g,~ \text{ for } x, y\in g.
\end{eqnarray*}
\hfill $\square$

A Leibniz algebra $(g, [\c, \c]_g)$ is said to be nondegenerate if the bracket $[\c, \c]_g$ satisfy the followings
\begin{eqnarray*}
&&[x, y]_g = 0, ~\text{ for all } y  \text{ implies that } x= 0,\\
&&[x,y]_g = 0, ~ \text{ for all } x \text{ implies that } y=0.
\end{eqnarray*}

\begin{cor}
Let $(g, [\c, \c]_g)$ be a nondegenerate Leibniz algebra.  Let $R: g\rightarrow g$ be a  left $g$-linear map and $S: g\rightarrow g$ be a  right $g$-linear map. Then $(R, S)$ is a Rota-Baxter system on $(g, [\c, \c]_g)$ if and only if
\begin{eqnarray*}
R\circ S=S\circ R=0.
\end{eqnarray*}
\end{cor}

\medskip

Another class of Rota-Baxter systems arise from twisted Rota-Baxter operators. Let $(g, [\c, \c]_g)$ be a Leibniz algebra and $\sigma: g\rightarrow g$ be a Leibniz algebra morphism.

\begin{definition}
A linear map $R: g\rightarrow g$ is said to be a $\sigma$-twisted Rota-Baxter operator if $R$ satisfies
\begin{eqnarray}
[Rx, Ry]_g=R([Rx, y]_g+[x, (\sigma\circ R)y]_g), ~\text{ for all } x, y\in g.
\end{eqnarray}
\end{definition}
When $\sigma=Id$, a $\sigma$-twisted Rota-Baxter operator is nothing but a Rota-Baxter operator.

\begin{example}
A differential Rota-Baxter Leibniz algebra of weight $\lambda$ is a Leibniz algebra $(g, [\c, \c]_g)$ together with linear maps $R, \partial: g\rightarrow g$ satisfying the following set of identities
\begin{eqnarray*}
(dR1)&&~~ [Rx, Ry]_g=R([Rx, y]_g+[x, Ry]_g+\lambda[x, y]_g),\\
(dR2)&&~~ \partial[x, y]_g=[\partial x, y]_g+[x, \partial y]_g+\lambda[\partial x, \partial y]_g,\\
(dR3)&&~~\partial \circ R=Id.
\end{eqnarray*}
Let $(g, R, \partial)$ be a differential Rota-Baxter Leibniz algebra of weight $\lambda$. It follows from $(dR2)$ that
the map
\begin{eqnarray*}
\sigma: g\rightarrow g,~\sigma(x)=x+\lambda\partial(x),~\text{ for }  x\in g
\end{eqnarray*}
is a Leibniz algebra morphism. On the other hand, $(dR3)$ implies that
\begin{eqnarray*}
(\sigma\circ R)(x)=R(x)+\lambda x,~ \text{ for } x\in g.
\end{eqnarray*}
Hence, by $(dR2)$, we get
\begin{eqnarray*}
[Rx, Ry]_g=R([Rx, y]_g+[x, (\sigma\circ R)y]_g), ~\text{ for } x, y\in g.
\end{eqnarray*}
This shows that $R$ is a $\sigma$-twisted Rota-Baxter operator.
\end{example}
\begin{proposition}
Let $R$ be a $\sigma$-twisted Rota-Baxter operator on a Leibniz algebra $(g, [\c, \c]_g)$. Then $(R, \sigma \circ R)$ is a Rota-Baxter system on $(g, [\c, \c]_g)$.
\end{proposition}
{\bf Proof.} Note that the condition Eq. (3.1) is same as the first condition of a Rota-Baxter system. To prove the second one, we observe that
\begin{eqnarray*}
[(\sigma\circ R)x, (\sigma\circ R)y]_g&=& \sigma [Rx, Ry]_g\\
&=& (\sigma\circ R)([Rx, y]_g+[x, (\sigma\circ R)y]_g).
\end{eqnarray*}
This shows that $(R, \sigma\circ R)$ is a Rota-Baxter system on $(g, [\c, \c]_g)$. \hfill $\square$
\begin{example}
Let $(W, [\c, \c]_W)$ be the Witt Lie algebra generated by basis elements $\{l_n\}_{n\in \mathbb{Z}}$ and  the Lie bracket given by
\begin{eqnarray*}
[l_m, l_n]_{W}=(m-n)l_{m+n},~ \text{ for } m, n\in \mathbb{Z}.
\end{eqnarray*}
View this Lie algebra as a Leibniz algebra. Let $q\in \mathbb{K}$ be a nonzero scalar that is not a root of unity. We define linear maps  $\sigma, R: W\rightarrow W$ by
\begin{eqnarray*}
\sigma(l_n)=q^{n}l_n, ~~~R(l_n)=\frac{1-q}{1-q^{n}}l_n, ~ \text{ for } n\in \mathbb{Z}.
\end{eqnarray*}
Then $\sigma$ is clearly a Leibniz algebra morphism. Moreover, it is easy to verify that $R$ satisfies
\begin{eqnarray*}
[R(l_m), R(l_n)]_W=R([R(l_m), l_n]_W+[l_m, (\sigma\circ R)(l_n)]_W), ~ \text{ for } m, n\in \mathbb{Z}.
\end{eqnarray*}
Therefore, $R$ is a $\sigma$-twisted Rota-Baxter operator. Hence, $(R, \sigma\circ R)$ is a Rota-Baxter system on $W$.
\end{example}

\medskip

In \cite{PO17} the authors introduced a notion of weak pseudotwistor on an associative algebra and showed that a weak pseudotwistor induce a new associative algebra structure. A Rota-Baxter system on an associative algebra gives rise to a weak pseudotwistor, hence a new associative algebra structure. This is not true for Rota-Baxter systems on Leibniz algebras. However, if we concentrate on Rota-Baxter operators, they induce a new Leibniz algebra structure via a Leibniz analogue of weak pseudotwistor. Let us first recall the new Leibniz algebra associated to a Rota-Baxter operator on a Leibniz algebra.

Let $(g, [\c, \c]_g)$ be a Leibniz algebra, and $R: g \rightarrow g$ be a Rota-Baxter operator, i.e., $R$ satisfies
$$[Rx, Ry]_g = R ( [Rx, y] + [x, Ry]), ~ \text{ for } x, y \in g.$$
Then the vector space $g$ carries a new Leibniz algebra structure with bracket $$[x, y]_R = [Rx, y] + [x, Ry], ~ \text{ for } x, y \in g.$$

Here we give a new example of Rota-Baxter operator on a Leibniz algebra induced from dialgebra \cite{L01}.

\begin{definition}
A dialgebra is a vector space $D$ together with two bilinear operations $\dashv, \vdash: D\o D\rightarrow D$ satisfying the following identities
\begin{eqnarray*}
&& a\dashv (b \dashv c)=(a\dashv b)\dashv c=a\dashv (b \vdash c),\\
&& (a\vdash b)\dashv c=a\vdash (b\dashv c),\\
&& (a\dashv b)\vdash c=(a\vdash b)\vdash c=a\vdash (b \vdash c),~ \text{ for } a, b, c\in D.
\end{eqnarray*}
\end{definition}
A dialgebra as above may be denoted by the triple $(D, \dashv, \vdash)$. Any associative algebra is obviously a dialgebra with both the bilinear maps coincide with the associative product. See Loday \cite{L01} for more examples of dialgebras.
\begin{proposition}
Let $(D, \dashv, \vdash)$ be a dialgebra. Then $(D, [\c, \c]_D)$ is a Leibniz algebra, where
\begin{eqnarray*}
[a, b]_D:=a\vdash b-b\dashv a,~ \text{ for } a, b\in D.
\end{eqnarray*}
\end{proposition}
{\bf Proof.} For any $a, b, c\in D$, we have
\begin{eqnarray*}
&& [[a, b]_D, c]_D+[b, [a, c]_D]_D\\
&=& [a\vdash b-b\dashv a, c]_D+[b, a\vdash c-c\dashv a]_D\\
&=& (a\vdash b-b\dashv a)\vdash c- c \dashv (a\vdash b-b\dashv a)+ b \vdash (a\vdash c-c\dashv a)- (a\vdash c-c\dashv a) \dashv b\\
&=& a\vdash (b\vdash c-c\dashv b)- (b\vdash c-c\dashv b) \dashv a\\
&=&[a, [b, c]_D]_D.
\end{eqnarray*}
Hence $(D, [\c, \c]_D)$ is a Leibniz algebra. \hfill $\square$

The Leibniz algebra in the above proposition is called the Leibniz algebra induced from the dialgebra $(D, \dashv, \vdash)$.

\begin{definition}
Let $(D, \dashv, \vdash)$ be a dialgebra.  A Rota-Baxter operator on $D$ consists of a linear map $R: D\rightarrow D$ satisfying
\begin{eqnarray*}
&& R(a)\ast R(b)=R(R(a)\ast b+a\ast R(b)),
\end{eqnarray*}
for all $a, b\in D$ and $\ast=~\dashv, \vdash$.
\end{definition}
\begin{proposition}
Let $(D, \dashv, \vdash)$ be a dialgebra and $R$ be a Rota-Baxter operator on it. Then $R$ is a Rota-Baxter operator on the induced Leibniz algebra $(D, [\c, \c]_D)$.
\end{proposition}
{\bf Proof.} For any $a, b\in D$, we have
\begin{eqnarray*}
[Ra, Rb]_D&=& Ra\vdash Rb-Rb\dashv Ra\\
&=& R(R(a)\vdash b+a\vdash R(b))-R(R(b)\dashv a+b\dashv R(a))\\
&=&R([Ra, b]_D+[a, Rb]_D).
\end{eqnarray*}
Hence the result follows.
              \hfill $\square$

\medskip

The Leibniz bracket $[\c, \c]_R$ induced from a Rota-Baxter operator $R$ can be understood in terms of the weak pseudotwistor on a Leibniz algebra.
\begin{definition}
Let $(g, [\c,\c]_g)$ be a Leibniz algebra with the Leibniz bracket denoted by the product $\mu$. A linear map $T: g\o g\rightarrow g\o g$ is said to be  a weak pseudotwistor if there exists a linear map $\tau: g\o g\o g\rightarrow g\o g\o g$ with $(\eta_{12}\o Id)\circ \tau=\tau\circ (\eta_{12}\o Id)$ and commuting the following diagram
$$\aligned
\xymatrix
{
& g\o g \o g \ar[rr]^-{Id \o \mu} & & g \o g \ar[dd]^-{T} && g\o g \o g \ar[ll]_-{\mu \o Id}  &\\
g\o g \o g \ar[ru]^-{Id \o T}\ar[rd]_-{\tau} & & & && & g\o g \o g \ar[ld]^-{\tau} \ar[lu]_-{T \o Id}\\
& g\o g \o g \ar[rr]_-{Id \o \mu} & & g \o g  && g\o g \o g \ar[ll]^-{\mu \o Id}  &
}
\endaligned$$
Here $\eta_{12}:g\o g\rightarrow g\o g$ is the flip map $\eta_{12}(x\o y)=y\o x$. The map $\tau$ is called a weak companion of $T$.
\end{definition}
\begin{proposition}
Let $(g, [\c,\c]_g)$ be a Leibniz algebra and $T: g\o g\rightarrow g\o g$  be  a weak pseudotwistor. Then $(g, \mu \circ T)$ is a new Leibniz algebra structure on $g$.
\end{proposition}
{\bf Proof.} We have
\begin{eqnarray*}
&& (\mu \circ T)\circ (Id\o (\mu \circ T))\\
&=& \mu\circ(Id\o \mu ) \circ \tau\\
&=& \mu\circ( \mu  \o Id ) \circ \tau + \mu\circ(Id\o \mu ) \circ (\eta_{12}\o Id)\circ \tau\\
&=& (\mu \circ T)\circ ((\mu \circ T)\o Id)+\mu\circ(Id\o \mu )\circ \tau \circ (\eta_{12}\o Id)\\
&=& (\mu \circ T)\circ ((\mu \circ T)\o Id)+(\mu \circ T)\circ (Id\o (\mu \circ T))\circ (\eta_{12}\o Id).
\end{eqnarray*}
This shows that $\mu\circ T$ defines a Leibniz bracket on $g$.\hfill $\square$

\begin{proposition}
Let $(g, [\c,\c]_g)$ be a Leibniz algebra and $R: g\rightarrow g$ be a Rota-Baxter operator on it. Then the map $T: g\o g\rightarrow g\o g$ defined by
\begin{eqnarray*}
T(x\o y)=R(x)\o y+x\o R(y)
\end{eqnarray*}
is a weak pseudotwistor on $g$. Consequently, $g$ carries a new Leibniz algebra structure with bracket $[x, y]_R = [Rx, y]_g + [x, Ry]_g$, for $x, y \in g$.
\end{proposition}
{\bf Proof.} We define $\tau: g\o g\o g\rightarrow g\o g\o g$ by
\begin{eqnarray*}
\tau (x \o y \o z) = R(x) \o R(y) \o z+R(x) \o y \o R(z)+x \o R(y) \o R(z), ~\text{ for } x, y, z \in g.
\end{eqnarray*}
We will show that $T$ is a weak pseudotwistor with a weak companion $\tau$. First observe that
\begin{eqnarray*}
&&((\eta_{12}\o Id)\circ \tau)(x\o y\o z)\\
&=& R(y)\o R(x)\o z+y\o R(x)\o R(z)+R(y)\o x\o R(z)\\
&=& \tau(y\o x\o z)=(\tau \circ(\eta_{12}\o Id))(x\o y\o z).
\end{eqnarray*}
Next, we have
\begin{eqnarray*}
&&(T\circ (Id\o \mu \circ T))(x\o y\o z)\\
&=& R(x)\o \mu(R(y)\o z+ y\o R(z))+x\o \mu(R(y)\o R(z))\\
&=& (( Id\o \mu)\circ \tau )(x\o y\o z).
\end{eqnarray*}
Similarly, we have
\begin{eqnarray*}
T\circ ((\mu \circ T) \o Id) =( \mu \o Id)\circ \tau.
\end{eqnarray*}
Hence, the result follows.\hfill $\square$

\begin{remark}
The notion of weak pseudotwistor on a Leibniz algebra is a generalization of weak pseudotwistor on an associative algebra introduced by  Panaite and  Oystaeyen \cite{PO17}.
In the associative context, a  Rota-Baxter system induces a weak pseudotwistor on the underlying associative algebra. It is remarked that given a Leibniz algebra $(g, [\c,\c]_g)$
and a Rota-Baxter system $(R, S)$  on $g$, the map
\begin{eqnarray*}
T: g\o g\rightarrow g\o g, ~~T(x\o y)=R(x)\o y+x\o S(y)
\end{eqnarray*}
is not a weak pseudotwistor on $g$ with weak companion
\begin{eqnarray*}
\tau(x\o y\o z)=R(x)\o R(y)\o z+R(x)\o y\o S(z)+x\o S(y)\o S(z)
\end{eqnarray*}
as $(\eta_{12}\o Id)\circ \tau\neq\tau\circ (\eta_{12}\o Id)$.
\end{remark}
\section{Maurer-Cartan characterization of relative Rota-Baxter systems}
\def\theequation{\arabic{section}. \arabic{equation}}
\setcounter{equation} {0}

In the section, we construct a graded Lie algebra that characterize relative Rota-Baxter systems as Maurer-Cartan elements. Using this characterization, we define the cohomology associated to a relative Rota-Baxter system. We first recall some results from \cite{B97}.

A permutation $\sigma\in \mathbb{S}_n$ is called an $(i, n-i)$-shuffle if $\sigma(1) < \c \c \c <\sigma(i)$ and $\sigma(i+ 1) <
\c \c \c < \sigma(n)$. If $i = 0$ or $n$ we assume $\sigma = Id$. The set of all $(i, n-i)$-shuffles will be
denoted by $\mathbb{S}_{(i,n-i)}$.

Let $M$ be a vector space. We consider the graded vector space
$$C^{\ast}(M, M)=\oplus_{n\geq1}C^n(M,M)=\oplus_{n\geq1}Hom(\otimes^{n}M, M)$$
of multilinear maps on $M$. The Balavoine bracket is a degree $-1$ bracket on the graded vector space $C^{\ast}(M, M)$ given by
\begin{eqnarray*}
[f, g]_B := f\overline{\circ}g- (-1)^{pq}g\overline{\circ}f,
\end{eqnarray*}
for $f\in C^{p+1}(M, M), g\in C^{q+1}(M, M).$ Here $f\overline{\circ}g\in  C^{p+q+1}(M, M)$ is defined by
\begin{eqnarray*}
f\overline{\circ}g=\sum_{k=1}^{p+1}(-1)^{(k-1)q}f\circ_k g,
\end{eqnarray*}
with
\begin{eqnarray*}
&&(f\circ_k g)(x_1,\c \c \c, x_{p+q+1})\\
&=& \sum_{\sigma\in \mathbb{S}_{(k-1, q)}}(-1)^{\sigma}f(x_{\sigma(1)}, \c \c \c, x_{\sigma(k-1)}, g(x_{\sigma(k)}, \c \c \c, x_{\sigma(k+q-1)}, x_{k+q}), x_{k+q+1}, \c \c \c, x_{p+q+1}).
\end{eqnarray*}

\begin{theorem}(\cite{B97})
With the above notations, $(C^{\ast}(M, M), [\c, \c]_B)$ is a degree $-1$ graded Lie algebra. In other words $(C^{\ast +1} (M,M), [\c, \c]_B)$ is a graded Lie algebra. Its
Maurer-Cartan elements are precisely the Leibniz algebra structures on $M$.

\end{theorem}

Let $(V, \rho^L, \rho^R)$ be a representation of a Leibniz algebra $(g, [\c, \c]_g)$. Consider the semidirect product
Leibniz algebra structure on $g \oplus g \oplus V$. We denote the corresponding Leibniz product by $\widehat{\mu}$. Then $\widehat{\mu}$ is a Maurer-Cartan element in the graded Lie algebra $(C^{\ast +1} (g \oplus g \oplus V, g \oplus g \oplus V), [\c, \c]_B)$.

Consider the graded vector subspace $C^\ast (V, g) \subset C^\ast ( g \oplus g \oplus V, g \oplus g \oplus V )$ given by
\begin{eqnarray*}
C^{\ast}(V, g):=\oplus_{n\geq 1}C^{n}(V, g):=\oplus_{n\geq 1}Hom(V^{\otimes n}, g\oplus g).
\end{eqnarray*}

\begin{theorem}
With the above notations, $(C^{\ast}(V, g), [[\c, \c]])$ is a graded Lie algebra, where the
graded Lie bracket $[[\c, \c]]: C^{m}(V, g)\times C^{n}(V, g)\rightarrow C^{m+n}(V, g)$ is defined by
\begin{eqnarray*}
[[(P, Q),(P', Q')]] := (-1)^{m}[[\widehat{\mu},(P, Q)]_B,(P', Q')]_B,
\end{eqnarray*}
for any $(P, Q)\in C^{m}(V, g), (P', Q')\in C^{n}(V, g)$. Moreover, its Maurer-Cartan elements are relative Rota-Baxter systems on the Leibniz algebra
$(g, [\c, \c]_g)$ with respect to the representation $(V, \rho^L, \rho^R)$.
\end{theorem}

 Let $Pr_1, Pr_2 : g\oplus g\rightarrow g$ denote the projection maps onto the first and second factor, respectively. Then the explicit description of the above graded Lie bracket is given by
\begin{small}
\begin{eqnarray*}
&&Pr_1([[(P, Q),(P', Q')]](v_1, \c \c \c, v_{m+n}))\\
&=& \sum_{k=1}^{m}\sum_{\sigma\in \mathbb{S}_{(k-1, n)}}(-1)^{(k-1)n}(-1)^{\sigma}P(v_{\sigma(1)}, \c \c \c, v_{\sigma(k-1)},\rho^{L}(P'(v_{\sigma(k)}, \c \c \c, v_{\sigma(k+n-1)}))v_{k+n},\c \c \c, v_{m+n})\\
&&+\sum_{k=2}^{m}\sum_{\sigma\in \mathbb{S}_{(k-2, n, 1)}}(-1)^{kn}(-1)^{\sigma}P(v_{\sigma(1)}, \c \c \c, v_{\sigma(k-2)},\rho^{R}(Q'(v_{\sigma(k)}, \c \c \c, v_{\sigma(k+n-2)}))\\
&&\hspace{10cm}v_{\sigma(k+n-1)}, v_{k+n},\c \c \c, v_{m+n})\\
&&+\sum_{k=1}^{n}\sum_{\sigma\in \mathbb{S}_{(k-1, m)}}(-1)^{(k+n-1)m}(-1)^{\sigma}P'(v_{\sigma(1)}, \c \c \c, v_{\sigma(k-1)},\rho^{L}(P(v_{\sigma(k)}, \c \c \c, v_{\sigma(k+m-1)}))\\
&&\hspace{11cm}v_{_{\sigma(k+m)}}, \c \c \c, v_{m+n})\\
&&\sum_{k=1}^{n}\sum_{\substack{\sigma\in \mathbb{S}_{(k-1, m,1)}, \\ \sigma(k+m-1)=k+m}}(-1)^{(k+n-1)m+1}(-1)^{\sigma}P'(v_{\sigma(1)}, \c \c \c, v_{\sigma(k-1)},\rho^{R}(Q(v_{\sigma(k)}, \c \c \c, v_{\sigma(k-1+m)}))\\
&&\hspace{10cm}v_{\sigma(k+m)}, v_{k+m+1}, \c \c \c, v_{m+n})\\
&&+ \sum_{\sigma\in \mathbb{S}_{(m,n-1)}}(-1)^{mn+1}(-1)^{\sigma}[P(v_{\sigma(1)}, \c \c \c, v_{\sigma(m)}), P'(v_{\sigma(m+1)}, \c \c \c, v_{\sigma(m+n-1)}, v_{m+n}]_g\\
&&+\sum_{k=1}^{m}\sum_{\sigma\in \mathbb{S}_{(k-1, n-1)}}(-1)^{(k-1)n}(-1)^{\sigma}[P'(v_{\sigma(k)}, \c \c \c, v_{\sigma(k+n-2)}), P(v_{\sigma(1)}, \c \c \c,v_{\sigma(k-1)}, v_{k+n}, \c \c \c, v_{m+n})]_g,
\end{eqnarray*}
\end{small}
and
\begin{small}
\begin{eqnarray*}
&&Pr_2([[(P, Q),(P', Q')]](v_1, \c \c \c, v_{m+n}))\\
&=& \sum_{k=1}^{m}\sum_{\sigma\in \mathbb{S}_{(k-1, n)}}(-1)^{(k-1)n}(-1)^{\sigma}Q(v_{\sigma(1)}, \c \c \c, v_{\sigma(k-1)},\rho^{L}(P'(v_{\sigma(k)}, \c \c \c, v_{\sigma(k+n-1)}))v_{k+n},\c \c \c, v_{m+n})\\
&&+\sum_{k=2}^{m}\sum_{\sigma\in \mathbb{S}_{(k-2, n, 1)}}(-1)^{kn}(-1)^{\sigma}Q(v_{\sigma(1)}, \c \c \c, v_{\sigma(k-2)},\rho^{R}(Q'(v_{\sigma(k)}, \c \c \c, v_{\sigma(k+n-2)}))\\
&&\hspace{10cm}v_{\sigma(k+n-1)}, v_{k+n},\c \c \c, v_{m+n})\\
&&+\sum_{k=1}^{n}\sum_{\sigma\in \mathbb{S}_{(k-1, m)}}(-1)^{(k+n-1)m}(-1)^{\sigma}Q'(v_{\sigma(1)}, \c \c \c, v_{\sigma(k-1)},\rho^{L}(P(v_{\sigma(k)}, \c \c \c, v_{\sigma(k+m-1)}))\\
&&\hspace{11cm}v_{_{\sigma(k+m)}}, \c \c \c, v_{m+n})\\
&&\sum_{k=1}^{n}\sum_{\substack{\sigma\in \mathbb{S}_{(k-1, m,1)}, \\ \sigma(k+m-1)=k+m}}(-1)^{(k+n-1)m+1}(-1)^{\sigma}Q'(v_{\sigma(1)}, \c \c \c, v_{\sigma(k-1)},\rho^{R}(Q(v_{\sigma(k)}, \c \c \c, v_{\sigma(k-1+m)}))\\
&&\hspace{10cm}v_{\sigma(k+m)}, v_{k+m+1}, \c \c \c, v_{m+n})\\
&&+ \sum_{\sigma\in \mathbb{S}_{(m,n-1)}}(-1)^{mn+1}(-1)^{\sigma}[Q(v_{\sigma(1)}, \c \c \c, v_{\sigma(m)}), Q'(v_{\sigma(m+1)}, \c \c \c, v_{\sigma(m+n-1)}, v_{m+n}]_g\\
&&+\sum_{k=1}^{m}\sum_{\sigma\in \mathbb{S}_{(k-1, n-1)}}(-1)^{(k-1)n}(-1)^{\sigma}[Q'(v_{\sigma(k)}, \c \c \c, v_{\sigma(k+n-2)}), Q(v_{\sigma(1)}, \c \c \c,v_{\sigma(k-1)}, v_{k+n}, \c \c \c, v_{m+n})]_g,
\end{eqnarray*}
\end{small}
for any $(P, Q)\in C^{m}(V, g), (P', Q')\in C^{n}(V, g)$.

{\bf Proof.} The graded Lie algebra  $(C^{\ast}(V, g), [[\c, \c]])$  is obtained via the derived bracket \cite{ST19}.   First consider the  graded
Lie algebra $(C^{\ast + 1}(g\oplus g\oplus V, g\oplus g\oplus V), [\c, \c]_B)$. Since $\widehat{\mu}$ is the semidirect product Leibniz algebra structure
on the vector space $g\oplus g\oplus V$, we deduce that $(C^{\ast + 1}(g\oplus g\oplus V, g\oplus g\oplus V), [\c, \c]_B, d=[\widehat{\mu},\c]_B)$  is a differential graded Lie algebra. Obviously $C^{\ast + 1}(V, g)$ is an abelian subalgebra. Therefore, by the derived bracket construction, we define a bracket on the shifted graded vector space $C^{\ast}(V, g)$ by
\begin{eqnarray*}
[[(P, Q),(P', Q']] :=(-1)^{m}[d((P, Q)), (P', Q')]_B= (-1)^{m}[[\widehat{\mu},(P, Q)],(P', Q')],
\end{eqnarray*}
for any $(P, Q)\in C^{m}(V, g), (P', Q')\in C^{n}(V, g)$.
The derived bracket $[[\c, \c]]$ is closed on $C^{\ast}(V, g)$, which implies that ($C^{\ast}(V, g), [[\c, \c]]$) is a graded Lie algebra.

\medskip

For $ (R, S)\in C^1(V, g)$, we have
\begin{eqnarray*}
&& Pr_1([[(R, S),(R, S)]](u, v))=2([Ru, Rv]_g-R(\rho^{L}(Ru)v)-R(\rho^{R}(Sv)u)),\\
&& Pr_2([[(R, S),(R, S)]](u, v))=2([Su, Sv]_g-S(\rho^{L}(Ru)v)-S(\rho^{R}(Sv)u)).
\end{eqnarray*}
\noindent Thus, $(R, S)$ is a  Maurer-Cartan element  (i.e. $[[(R, S),(R, S)]] = 0$) if and only if $(R, S)$ is a  relative Rota-Baxter systems on $g$ with respect to the representation $(V, \rho^L, \rho^R)$. The proof is finished.  \hfill $\square$

Thus, relative Rota-Baxter systems can be characterized as Maurer-Cartan elements in a graded Lie algebra. It
follows from the above theorem that if $(R, S)$ is a relative Rota-Baxter system, then $d_{(R,S)}:= [[(R, S), \c]]$
is a differential on $C^{\ast}(V, g)$ and makes the gLa $(C^{\bullet}(V, g), [[\c , \c]])$ into a differential graded Lie algebra.

\medskip

The cohomology of the cochain complex $(C^{\bullet}(V, g), d_{(R,S)})$ is called the cohomology of the relative
Rota-Baxter system $(R, S)$. We denote the corresponding cohomology groups simply by $H^{\bullet}(V, g)$.

The following theorem describes the Maurer-Cartan deformation of a relative Rota-Baxter system.

\begin{theorem}
Let $(R, S)$ be a relative Rota-Baxter system on a Leibniz algebra $(g, [\c, \c]_g)$ with respect to a representation $(V, \rho^L , \rho^{R})$.  For any pair $(R', S')$ of linear maps from $V$ to $g$, the pair of sums $(R + R', S + S')$ is a relative Rota-Baxter system if and only if $(R', S')$ is a Maurer-Cartan element in the differential graded Lie algebra $(C^{\ast}(V, g), [[\c, \c]], d_{(R,S)})$, i.e.
\begin{eqnarray*}
[[(R + R', S + S'),(R + R', S + S')]]=0 \Leftrightarrow d_{(R,S)}(R', S')+\frac{1}{2}[[(R', S'),(R', S')]]=0.
\end{eqnarray*}
\end{theorem}
\section{Deformations of relative Rota-Baxter systems }
\def\theequation{\arabic{section}. \arabic{equation}}
\setcounter{equation} {0}

\subsection{Formal deformations}

Let $\mathbb{K}[[t]]$ be the ring of power series in one variable $t$. For any $\mathbb{K}$-linear space $V$,
let $V [[t]]$ denotes the vector space of formal power series in $t$ with coefficients from $V$.
If in addition, $(g, [\c, \c]_g)$ is a Leibniz algebra over $\mathbb{K}$, then there is a $\mathbb{K}[[t]]$-Leibniz
algebra structure on $g[[t]]$ given by
\begin{eqnarray*}
[\sum_{i=0}^{+\infty}x_it^{i},\sum_{j=0}^{+\infty}y_jt^{j} ]_g=\sum_{k=0}^{+\infty}\sum_{i+j=k}[x_i, y_j]t^{k}, ~ \text{ for all } x_i, y_j\in g.
\end{eqnarray*}
Let $(V, \rho^L, \rho^{R})$ be a representation of the Leibniz algebra $(g, [\c, \c]_g)$. Then there is a representation $(V [[t]], \rho^L, \rho^{R})$ of the $\mathbb{K}[[t]]$-Leibniz algebra $g[[t]]$. Here $\rho^L$ and $\rho^R$ are given by
\begin{eqnarray*}
&&\rho^L(\sum_{i=0}^{+\infty}x_it^{i})(\sum_{j=0}^{+\infty}v_jt^{j})=\sum_{k=0}^{+\infty}\sum_{i+j=k}\rho^L(x_i)(v_j)t^{k},\\
&&\rho^R(\sum_{i=0}^{+\infty}x_it^{i})(\sum_{j=0}^{+\infty}v_jt^{j})=\sum_{k=0}^{+\infty}\sum_{i+j=k}\rho^R(x_i)(v_j)t^{k}, ~ \text{ for all } x_i\in g, v_j\in V.
\end{eqnarray*}
Let $(R, S)$ be a relative Rota-Baxter system on the Leibniz algebra $(g, [\c, \c]_g)$ with respect to the representation $(V, \rho^L, \rho^{R})$.  We consider two power series
\begin{eqnarray*}
R_t=\sum_{i=0}^{+\infty}\mathfrak{R}_it^{i} \text{ and } S_t=\sum_{j=0}^{+\infty}\mathfrak{S}_jt^{j}, \text{ where }
\mathfrak{R}_i, \mathfrak{S}_j \in Hom_{\mathbb{K}}(V, g).
\end{eqnarray*}
That is, both $R_t$ and  $S_t$ are in  $Hom_{\mathbb{K}}(V, g)[[t]]$. Extend them to  $\mathbb{K}[[t]]$-linear maps from $V [[t]]$ to
$g[[t]]$. We still denote them by same symbols.

\begin{definition}
If $R_t=\sum_{i=0}^{+\infty}\mathfrak{R}_it^{i}$ and  $S_t=\sum_{j=0}^{+\infty}\mathfrak{S}_jt^{j}$ with $\mathfrak{R}_0 = R$, $\mathfrak{S}_0 = S$ satisfy
\begin{eqnarray*}
&&[R_tu, R_tv]_g=R_t(\rho^{L}(R_tu)v+\rho^{R}(S_tv)u),\\
&&[S_tu, S_tv]_g=S_t(\rho^{L}(R_tu)v+\rho^{R}(S_tv)u)),
\end{eqnarray*}
we say that $(R_t, S_t)$ is a formal deformation of the relative Rota-Baxter system $(R, S)$.
\end{definition}

By expanding these equations and comparing coefficients of various powers of $t$, we obtain for $k\geq 0$,
\begin{eqnarray*}
 &&\sum_{k=0}^{+\infty}\sum_{i+j=k}[\mathfrak{R}_iu, \mathfrak{R}_jv]_g=\sum_{k=0}^{+\infty}\sum_{i+j=k}\mathfrak{R}_i(\rho^{L}(\mathfrak{R}_ju)v+\rho^{R}(\mathfrak{S}_jv)u),\\
&&\sum_{k=0}^{+\infty}\sum_{i+j=k}[\mathfrak{S}_iu, \mathfrak{S}_jv]_g=\sum_{k=0}^{+\infty}\sum_{i+j=k}\mathfrak{S}_i(\rho^{L}(\mathfrak{R}_ju)v+\rho^{R}(\mathfrak{S}_jv)u)).
\end{eqnarray*}
Both of these identities hold for $k = 0$ as $(R, S)$ is a relative Rota-Baxter system. For $k = 1$, we get
\begin{eqnarray*}
 &&[Ru,  \mathfrak{R}_1v]_g +[\mathfrak{R}_1u,  Rv]_g =\mathfrak{R}_1(\rho^{L}(Ru)v+\rho^{R}(Sv)u)+R(\rho^{L}(\mathfrak{R}_1u)v+\rho^{R}(\mathfrak{S}_1v)u),\\
&&[Su,  \mathfrak{S}_1v]_g +[\mathcal{S}_1u,  Sv]_g =\mathfrak{S}_1(\rho^{L}(Ru)v+\rho^{R}(Sv)u)+S(\rho^{L}(\mathfrak{R}_1u)v+\rho^{R}(\mathfrak{S}_1v)u,
\end{eqnarray*}
for $u, v \in V$. These identities are equivalent to the single condition
\begin{eqnarray*}
[[(R, S), (\mathfrak{R}_1, \mathfrak{S}_1)]]=0.
\end{eqnarray*}

As a consequence, we get the following.

\begin{proposition}
Let $(R_t=\sum_{i=0}^{+\infty}\mathfrak{R}_it^{i}, S_t=\sum_{j=0}^{+\infty}\mathfrak{S}_jt^{j})$ be a formal deformation of a relative Rota-Baxter system $(R, S)$ on the Leibniz algebra $(g, [\c, \c]_g)$  with respect to a representation $(V, \rho^L, \rho^{R})$. Then $(\mathfrak{R}_1, \mathfrak{S}_1)$ is a $1$-cocycle in the cohomology of the relative Rota-Baxter system
$(R, S)$, that is, $d_{(R, S)}(\mathfrak{R}_1, \mathfrak{S}_1)=0$.
\end{proposition}

\begin{definition}
Let $(R, S)$ be a relative Rota-Baxter system  on the Leibniz algebra $(g, [\c, \c]_g)$  with respect to a representation $(V, \rho^L, \rho^{R})$. The 1-cocycle $(\mathfrak{R}_1, \mathfrak{S}_1)$ is called the infinitesimal of the formal deformation   $(R_t=\sum_{i=0}^{+\infty}\mathfrak{R}_it^{i}, S_t=\sum_{j=0}^{+\infty}\mathfrak{S}_jt^{j})$ of the relative Rota-Baxter system $(R, S)$.
\end{definition}
\begin{definition}
Two formal deformations $(R_t, S_t)$ and $(R'_t, S'_t)$  of a relative Rota-Baxter system $(R, S)$ on the Leibniz algebra $(g, [\c, \c]_g)$  with respect to a representation $(V, \rho^L, \rho^{R})$  are said to be equivalent if there exist two elements $x, y \in g$ and linear maps
$\phi_i, \varphi_i \in gl(g)$ and $\psi_i\in gl(V)$ for $i \geq 2$ such that for
\begin{eqnarray*}
&& \phi_t = Id_g + t(L_x-R_x) + \sum_{i=2}^{+\infty}\phi_i t^{i},~~~~~~ \varphi_t = Id_g + t(L_y-R_y) + \sum_{i=2}^{+\infty}\varphi_i t^{i}\\
&& \text{ and } \quad \psi_t=Id_V + t(\rho^{L}(x)-\rho^{R}(y)) + \sum_{i=2}^{+\infty}\psi_i t^{i},
\end{eqnarray*}
the following conditions hold:
\begin{eqnarray*}
&&(i)~~ [\phi_t(z), \phi_t(w)]_g=\phi_t([z, w]_g),~~~[\varphi_t(z), \varphi_t(w)]_g=\varphi_t([z, w]_g);\\
&&(ii)~~\psi_t (\rho^{L}(z)u)=\rho^{L}(\phi_t(z))\psi_t(u);\\
&&(iii)~~\psi_t (\rho^{R}(z)u) =\rho^{R}(\varphi_t(z))\psi_t(u);\\
&&(iv)~~R'_t\circ \psi_t (u)= \phi_t \circ R_t (u), ~~~S'_t\circ \psi_t (u)= \varphi_t \circ S_t (u),
\end{eqnarray*}
for all $z, w\in g$ and $u \in V$.
\end{definition}
By expanding the identities in (iv) and equating coefficients of $t$ from both sides, we obtain
\begin{eqnarray*}
(\mathfrak{R}_1, \mathfrak{S}_1)(u)-(\mathfrak{R}'_1, \mathfrak{S}'_1)(u)
&=& [R(u), x]_g-R(\rho^{R}(y)u)-[x, R(u)]_g+R(\rho^{L}(x)u)\\
&&+[S(u), y]_g-S(\rho^{R}(y)u)-[y, S(u)]_g+S(\rho^{L}(x)u)\\
&=&(d_{(R, S)}(x,y))(u).
\end{eqnarray*}
Thus, we have the following.
\begin{theorem}
The cohomology class of the infinitesimal of a deformation of a relative Rota-Baxter system depends only on the equivalence class of the deformation.
\end{theorem}


\subsection{Finite order deformations of a relative Rota-Baxter system}

In this subsection, we introduce a cohomology class associated to any order $n$ deformation of a relative
Rota-Baxter system, and show that an order $n$ deformation is extensible if and only if this cohomology class is trivial. Thus, we call this cohomology class the obstruction class of the order $n$ deformation being extensible.
\begin{definition}
Let $(R, S)$ be a relative Rota-Baxter system on a Leibniz algebra $(g, [\c, \c]_g)$ with respect to a representation $(V, \rho^L, \rho^{R})$. If the finite sums $$R_t=\sum_{i=0}^{n}\mathfrak{R}_it^{i} \text{ and }  S_t=\sum_{j=0}^{n}\mathfrak{S}_jt^{j} \text{ with } \mathfrak{R}_0 = R,~ \mathfrak{S}_0 = S$$ as $\mathbb{K}[[t]]/(t^{n+1})$-module maps from $V[[t]]/(t^{n+1})$ to the Leibniz algebra $g[[t]]/(t^{n+1})$ satisfy
\begin{eqnarray*}
&&[R_tu, R_tv]_g=R_t(\rho^{L}(R_tu)v+\rho^{R}(S_tv)u),\\
&&[S_tu, S_tv]_g=S_t(\rho^{L}(R_tu)v+\rho^{R}(S_tv)u)), ~\text{ for } u, v\in V,
\end{eqnarray*}
we say that $(R_t, S_t)$ is an order $n$ deformation of the relative Rota-Baxter system $(R, S)$.
\end{definition}

\begin{definition}
Let $(R_t, S_t)$ be an order $n$ deformation of the relative Rota-Baxter system $(R, S)$ on a Leibniz algebra $(g, [\c, \c]_g)$ with respect to a representation $(V, \rho^L, \rho^{R})$. If  there exists a pair $(\mathfrak{R}_{n+1}, \mathfrak{S}_{n+1})$ of linear maps from $V$ to $g$ such that $$(\widehat{R}_t = R_t +t^{n+1}\mathfrak{R}_{n+1}, \widehat{S}_t = S_t +t^{n+1}\mathfrak{S}_{n+1})$$ is a deformation of order $n+ 1$, we say that $(R_t, S_t)$ is extensible.
\end{definition}

Let $(R_t, S_t)$ be an order $n$ deformation of the relative Rota-Baxter system $(R, S)$ on a Leibniz algebra $(g, [\c, \c]_g)$ with respect to a representation $(V, \rho^L, \rho^{R})$. Define an element $Ob_{(R_t, S_t)}\in C^{2}(V, g)$  by
\begin{eqnarray}
Ob_{(R_t, S_t)}=-\frac{1}{2}\sum_{i+j=n+1, i, j\geq 1}[[(\mathfrak{R}_i, \mathfrak{S}_i),(\mathfrak{R}_j , \mathfrak{S}_j)]].
\end{eqnarray}
\begin{proposition}
The $2$-cochain $Ob_{(R_t, S_t)}$ is a 2-cocycle, that is, $d_{(R, S)}(Ob_{(R_t, S_t)})= 0$.
\end{proposition}
{\bf Proof.} We have
\begin{eqnarray*}
&&d_{(R, S)}(Ob_{(R_t, S_t)})\\
&=&-\frac{1}{2}\sum_{i+j=n+1, i, j\geq 1}[[(R, S), [[(\mathfrak{R}_i, \mathfrak{S}_i),(\mathfrak{R}_j , \mathfrak{S}_j)]]]]\\
&=&-\frac{1}{2}\sum_{i+j=n+1, i, j\geq 1}([[[[(R, S), (\mathfrak{R}_i, \mathfrak{S}_i)]],(\mathfrak{R}_j , \mathfrak{S}_j)]]-[[(\mathfrak{R}_i, \mathfrak{S}_i), [[(R, S),(\mathfrak{R}_j , \mathfrak{S}_j)]]]])\\
&=&\frac{1}{4}\sum_{i_1+i_2+j=n, i_1,i_2, j\geq 1}[[[[(\mathfrak{R}_{i_1}, \mathfrak{S}_{i_1}), (\mathfrak{R}_{i_2}, \mathfrak{S}_{i_2})]],(\mathfrak{R}_j , \mathfrak{S}_j)]]\\
&&-\frac{1}{4}\sum_{i+j_1+j_2=n, i,j_1, j_2\geq 1}[[(\mathfrak{R}_i, \mathfrak{S}_i), [[(\mathfrak{R}_{j_1}, \mathfrak{S}_{j_1}),(\mathfrak{R}_{j_2}, \mathfrak{S}_{j_2})]]]]\\
&=&\frac{1}{2}\sum_{i+j+k=n+1, i, j,k\geq 1}[[[[(\mathfrak{R}_i, \mathfrak{S}_i), (\mathfrak{R}_j, \mathfrak{S}_j)]], (\mathfrak{R}_k, \mathfrak{S}_k)]]\\
&=&0.
\end{eqnarray*}
The proof is finished. \hfill $\square$

\begin{definition}
Let $(R_t, S_t)$ be an order $n$ deformation of the relative Rota-Baxter system $(R, S)$ on a Leibniz algebra $(g, [\c, \c]_g)$ with respect to a representation $(V, \rho^L, \rho^{R})$.  The cohomology class $[Ob_{(R_t, S_t)}]\in H^{2}(V, g)$ is called the obstruction class for $(R_t, S_t)$ being extensible.
\end{definition}
As a consequence of Eq. (5.1) and Proposition 5.8, we obtain the following.
\begin{theorem}
Let $(R_t, S_t)$ be an order $n$ deformation of the relative Rota-Baxter system $(R, S)$ on a Leibniz algebra $(g, [\c, \c]_g)$ with respect to a representation $(V, \rho^L, \rho^{R})$. Then $(R_t, S_t)$ is extensible if and only if the obstruction class $[Ob_{(R_t, S_t)}]$ is trivial.
\end{theorem}
\begin{cor}
If $H^{2}(V, g)=0$ then every $1$-cocycle in the cohomology of a relative Rota-Baxter system $(R, S)$ is the infinitesimal of some formal deformation of $(R, S)$.
\end{cor}

\begin{center}
 {\bf ACKNOWLEDGEMENT}
 \end{center}

The work of A. Das is supported by the fellowship of Indian Institute of Technology (IIT) Kanpur. The work of S. Guo is supported by the NSF of China (No. 11761017) and Guizhou Provincial  Science and Technology  Foundation (No. [2020]1Y005).

\renewcommand{\refname}{REFERENCES}

\end{document}